\documentclass{article}\begin{document}
\centerline{\Large\bf On the Gauss Circle Problem}
\[
\]
\centerline{\bf Nikos Bagis}
\[
\]
\centerline{\bf Abstract}
We analyze the double series of Bessel functions given by Ramanujan. Using a very simple lemma we establish the uniform convergence of these series. By this we address to the Gauss circle problem.      
\[
\]        
\textbf{keywords}: Asymptotics, Summation, Ramanujan, Series, Integrals, Convergence  

\[
\]

\section{Introduction}
 
The general Bessel function of the first kind and $\nu$-th order is defined by
\begin{equation}
J_{\nu}(z):=\sum^{\infty}_{k=0}\frac{(-1)^k}{k!\Gamma(k+\nu+1)}\left(\frac{z}{2}\right)^{\nu+2k}\textrm{, }0\leq |z|<\infty\textrm{, }\nu\in \textbf{C}
\end{equation}
From this we define $Y_{\nu}(z)$ as
\begin{equation}
Y_{\nu}(z):=\frac{J_{\nu}(z)\cos(\nu\pi)-J_{-\nu}(z)}{\sin(\nu\pi)}
\end{equation}
and the modified Bessel function $K_{\nu}(z)$ by
\begin{equation}
K_{\nu}(z):=\frac{\pi}{2}\frac{e^{\pi i\nu/2}J_{-\nu}(iz)-e^{-\pi i\nu/2}J_{\nu}(iz)}{\sin(\nu\pi)}.
\end{equation}
If $\nu$ is an integer $n$, then it is understood that we define the functions above by taking the limits as $\nu\rightarrow n$. We also define $I_{n}(z)$ as 
\begin{equation}
I_{\nu}(z)=-Y_{\nu}(z)-\frac{2}{\pi}K_{\nu}(z)
\end{equation}
In particular we shall need only the functions $J_1(x)$ and $I_1(x)$.\\   
Following now Ramanujan we define $F(x)$ as
\begin{equation}
F(x)=
\left\{\begin{array}{cc}
 [x],   \textrm{ if }x\textrm{ is not an integer }\\
x-\frac{1}{2}, \textrm{ if }x \textrm{ is an integer }	
\end{array}\right\}
\end{equation}
where $[x]$ is the greatest integer less than or equal to $x$.\\
Then Ramanujan has claimed that (see [1]):\\
\\
\textbf{Conjecture 1.}\\
Let $F(x)$ be as defined by (5). If $\theta\in(0,1)$ and $x>0$, then
$$
\sum^{\infty}_{n=1}F\left(\frac{x}{n}\right)\sin(2\pi n \theta)=\pi x\left(\frac{1}{2}-\theta\right)-\frac{1}{4}\cot(\pi\theta)+
$$
\begin{equation}
+\frac{1}{2}\sqrt{x}\sum^{\infty}_{m=1}\sum^{\infty}_{n=0}\left(\frac{J_1(4\pi\sqrt{m(n+\theta)x})}{\sqrt{m(n+\theta)}}-\frac{J_1(4\pi\sqrt{m(n+1-\theta)x})}{\sqrt{m(n+1-\theta)}}\right)
\end{equation}
\\
\textbf{Conjecture 2.}\\
Let $F(x)$ be as defined by (5). If $\theta\in(0,1)$ and $x>0$, then
$$
\sum^{\infty}_{n=1}F\left(\frac{x}{n}\right)\cos(2\pi n \theta)=\frac{1}{4}-x\log(2\sin(\pi\theta))+
$$
\begin{equation}
+\frac{1}{2}\sqrt{x}\sum^{\infty}_{m=1}\sum^{\infty}_{n=0}\left(\frac{I_1(4\pi\sqrt{m(n+\theta)x})}{\sqrt{m(n+\theta)}}+\frac{I_1(4\pi\sqrt{m(n+1-\theta)x})}{\sqrt{m(n+1-\theta)}}\right)
\end{equation}
\\
These formulas  are connected with Gauss circle problem and Dirichlet divisors problem, that is, to determinate  respectively the best possible error terms $P(x)$ and $\Delta(x)$ in the following asymptotic expansions
\begin{equation}
\sum^{*}_{0\leq n\leq x}r_2(n)=\pi x+P(x)\textrm{, as }x\rightarrow\infty
\end{equation}                             
and
\begin{equation}
\sum^{*}_{n\leq x}d(n)=x\log(x)+(2\gamma-1)x+\frac{1}{4}+\Delta(x)\textrm{, as }x\rightarrow\infty.
\end{equation}
The asterisk $*$ in the summations means that if $x$ is an integer, only $\frac{1}{2}r_2(x)$, respectively $\frac{1}{2}d(x)$ is counted. Moreover $r_2(0)=1$.\\ 
Using (6) Berndt and Zaharescu [1] derived and prove the following representation
$$
\sum^{*}_{0\leq n\leq x}r_2(n)=\pi x+
$$
\begin{equation}
+2\sqrt{x}\sum^{\infty}_{n=0}\sum^{\infty}_{m=1}\left(\frac{J_1\left(4\pi\sqrt{m(n+\frac{1}{4})x}\right)}{\sqrt{m(n+\frac{1}{4})}}-\frac{J_1\left(4\pi\sqrt{m(n+\frac{3}{4})x}\right)}{\sqrt{m(n+\frac{3}{4})}}\right)
\end{equation} 
which shall investigate it.\\  
Actually Hardy and Voronoi (1904) have established the following formulas
\begin{equation}
\sum^{*}_{0\leq n\leq x}r_2(n)=\pi x+\sum^{\infty}_{n=1}r_2(n)\sqrt{\frac{x}{n}}J_1(2\pi\sqrt{nx})\textrm{, }x\geq0
\end{equation}
and
\begin{equation}
\sum^{*}_{n\leq x}d(n)=x\log(x)+x(2\gamma-1)+\frac{1}{4}+\sum^{\infty}_{n=1}d(n)\sqrt{\frac{x}{n}}I_1(4\pi\sqrt{nx})\textrm{, }x>0
\end{equation}
respectively

\section{Some Results}

We will need the following Lemmas\\
\\
\textbf{Lemma 1.} (see [2])\\
Let $f(x)$ be a function with Taylor series in $(-a,a)$, $a\geq 1$. Let also its Taylor series converges absolutely in $1$. Then exists a constant $c=c(f)$ such that
\begin{equation}
\sum^{M}_{k=1}f\left(\frac{1}{k}\right)=\int^{M}_{1}f\left(\frac{1}{t}\right)dt+c(f)+O\left(\frac{1}{M}\right)
\end{equation}    
Moreover 
\begin{equation}
c(f)=f(0)+f'(0)\gamma+\sum^{\infty}_{s=2}\frac{f^{(s)}(0)}{s!}\left(\zeta(s)-\frac{1}{s-1}\right).
\end{equation}
where $\zeta(s)$ is Riemann's zeta function.\\
\\
\textbf{Proof.}\\
See [2].\\
\\
\textbf{Note.}\\ We say that $G(x)=O_1\left(H(x)\right)$ if exist constants $C_1,C_2$ such that $C_1 C_2>0$ and
\begin{equation}
C_1H(x)\leq G(x)\leq C_2 H(x)\textrm{, for all }x\in[x_0,\infty)\textrm{, }x_0>0
\end{equation}
\\
\textbf{Lemma 2.}\\
If $f$ is such Lemma 1 and have derivatives $f^{(s)}(0)$ of the same sign (always non positive or non negative), then
\begin{equation}
c(f,x)=f(x)+O_1\left(\int^{1}_{0}f(xt)dt\right)
\end{equation}
and
\begin{equation}
c(f,x)=f(x)-\int^{1}_{0}f(xt)dt+O_1\left(\frac{1}{x^2}\int^{x}_{0}\int^{w}_{0}f(t)dtdw\right)
\end{equation} 
\textbf{Proof.}\\
Using the next inequalities
$$
\frac{1}{s+1}\leq-\zeta(s)+\frac{1}{s-1}+1\leq\frac{2}{s+1}
$$ 
\begin{equation}
\frac{-1}{s+1}+\frac{-3}{(s+1)(s+2)}\leq 
\zeta(s)-\frac{1}{s-1}-1\leq \frac{-1}{s+1}+\frac{-1/4}{(s+1)(s+2)}
\end{equation}
for $s=1,2,\ldots,$ in (14), we get the relations.\\
\\
\textbf{Lemma 3.} (see [1])\\
If $f(x)$ can be represented as a Fourier integral, $f(x)$ tends to 0 as $x\rightarrow\infty$, and $xf'(x)\in L^{p}(0,\infty)$ for some $p$ with $1<p\leq2$, then
\begin{equation}
\lim_{M\rightarrow\infty}\left(\sum^{M}_{n=1}f(n)-\int^{M}_{0}f(t)dt\right)=\lim_{M\rightarrow\infty}\left(\sum^{M}_{n=1}g(n)-\int^{M}_{0}g(t)dt\right)
\end{equation}  
where 
\begin{equation}
g(x):=2\int^{\infty}_{0}f(t)\cos(2\pi x t)dt.
\end{equation}
\\
\textbf{Corollary.} (see [1])
$$
\sum^{\infty}_{m=1}\frac{I_1(4\pi\sqrt{m(n+\theta)x})}{\sqrt{m(n+\theta)}}=
$$
\begin{equation}
=\frac{1}{\pi(n+\theta)\sqrt{x}}\lim_{M\rightarrow\infty}\left(\sum^{M}_{m=1}\sin\left(\frac{2\pi(n+\theta)x}{m}\right)-\int^{M}_{0}\sin\left(\frac{2\pi(n+\theta)x}{t}\right)dt\right)
\end{equation}
\\
\textbf{Theorem 1.}
$$
\sum^{\infty}_{n=0}\sum^{\infty}_{m=1}\frac{I_1(4\pi\sqrt{m(n+\theta)x})}{\sqrt{m(n+\theta)}}=
$$
\begin{equation}
=\sum^{\infty}_{n=0}\frac{\sin(2\pi(n+\theta)x)}{\pi(n+\theta)\sqrt{x}}
+O_1\left(\sum^{\infty}_{n=0}\frac{x^{-3/2}}{(n+\theta)^2}\right)=O\left(x^{-1/2}\right)
\end{equation}
\\
\textbf{Proof.}\\
Set $f(x)=\sin(x)=\sum^{\infty}_{n=0}\frac{(-1)^nx^{2n+1}}{(2n+1)!}=\sum^{\infty}_{n=0}a_n(x)$ and write $f_{+}(x)=\sum^{\infty}_{n=0}a_{2n}(x)$, $f_{-}(x)=\sum^{\infty}_{n=0}a_{2n+1}(x)$, we have two functions with derivatives of the same sign and $f(x)=f_{+}(x)+f_{-}(x)$. From Corollary and Lemmas 1,2 we get easily
$$
\sum^{\infty}_{m=1}\frac{I_1(4\pi\sqrt{m(n+\theta)x})}{\sqrt{m(n+\theta)}}=
$$
$$
=\frac{1}{\pi(n+\theta)\sqrt{x}}\lim_{M\rightarrow\infty}\left(\sum^{M}_{m=1}\sin\left(\frac{2\pi(n+\theta)x}{m}\right)-\int^{M}_{0}\sin\left(\frac{2\pi(n+\theta)x}{t}\right)dt\right)=
$$
$$
=\frac{1}{\pi(n+\theta)\sqrt{x}}[-\int^{1}_{0}\sin\left(\frac{2\pi(n+\theta)x}{t}\right)dt+\sin\left(2\pi(n+\theta)x\right)-
$$
$$
-\frac{1}{2\pi(n+\theta)}\int^{2\pi(n+\theta)}_{0}\sin(t)dt+O_1\left(\frac{1}{(2\pi(n+\theta))^2}\int^{2\pi(n+\theta)}_{0}\int^{w}_{0}\sin(t)dtdw\right)]=
$$
$$
=\frac{\sin(Y_n)}{\pi(n+\theta)\sqrt{x}}-\frac{1}{\pi(n+\theta)\sqrt{x}Y_n}
+O_1\left(\frac{1}{Y_n\pi(n+\theta)\sqrt{x}}\right)
$$
Were we make use of
$$
\int^{1}_{0}\sin\left(\frac{Y_n}{t}\right)dt=\sin(Y_n)-\textrm{Ci}(Y_n)Y_n
$$
with $Y_n=2\pi(n+\theta)x$ and
$$
\textrm{Ci}(x)=-\int^{\infty}_{x}\frac{\cos(t)}{t}dt
$$
is the cosine integral and the asymptotic formula 
\begin{equation}
\textrm{Ci}(x)\sim\frac{\sin(x)}{x}\left(\sum^{\infty}_{n=0}\frac{(-1)^n(2n)!}{x^{2n}}\right)-\frac{\cos(x)}{x}\left(\sum^{\infty}_{n=0}\frac{(-1)^n(2n+1)!}{x^{2n+1}}\right)
\end{equation}
as $x\rightarrow\infty$.

\section{Some Evaluations and Directions}

With the methods described in [1], we have from the following cosine transform
\begin{equation}
\int^{\infty}_{0}\frac{J_1\left(4\pi\sqrt{t(n+\theta)x}\right)}{\sqrt{t(n+\theta)}}\cos(2\pi t w)dt=\frac{1}{\pi(n+\theta)\sqrt{x}}\sin^2\left(\frac{\pi(n+\theta)x}{w}\right)
\end{equation} 
Hence in view of Lemma 3 the next formula is valid
$$
S_1(n)=\sum^{\infty}_{m=1}\frac{J_1\left(4\pi\sqrt{m(n+\theta)x}\right)}{\sqrt{m(n+\theta)}}-\frac{1}{2\pi(n+\theta)x^{1/2}}=
$$
$$
\frac{2}{\pi(n+\theta)x^{1/2}}\lim_{M\rightarrow\infty}\left(\sum^{M}_{k=1}\sin^2\left(\frac{\pi(n+\theta)x}{k}\right)-\int^{M}_{0}\sin^2\left(\frac{\pi(n+\theta)x}{t}\right)dt\right)=
$$
\begin{equation}
=\frac{2}{\pi(n+\theta)x^{1/2}}S(n,\theta)
\end{equation}
\\
But $\sin^2(x)=\frac{1-\cos(2x)}{2}$. Hence 
\begin{equation}
c(\sin^2,Y)=\sum^{\infty}_{n=1}\frac{(-1)^{n-1}2^{2n-1}Y^{2n}}{(2n)!}\zeta(2n)-\sum^{\infty}_{n=1}\frac{(-1)^{n-1}2^{2n-1}Y^{2n}}{(2n)!(2n-1)}
\end{equation}
and from [3] 
$$
\zeta(2s)=\frac{1}{\Gamma(2s)}\int^{\infty}_{0}\frac{t^{2s}}{t(e^t-1)}dt\textrm{, }Re(2s)>1
$$
Hence
$$
c(\sin^2,Y)=\int^{\infty}_{0}\sum^{\infty}_{n=1}\frac{(-1)^{n-1}2^{2n-1}Y^{2n}t^{2n}}{(2n)!\cdot\Gamma(2n)t(e^t-1)}dt+\sin^2(Y)-Y\textrm{Si}(2Y)=
$$ 
$$
=-\frac{\sqrt{Y}}{\sqrt{2}}\int^{\infty}_{0}\frac{\textrm{ber}^{(1)}(2\sqrt{2tY})}{\sqrt{t}(e^t-1)}dt+\sin^2(Y)-Y\textrm{Si}(2Y)
$$
Note that $\textrm{Si}(x)$ is the sin integral i.e
\begin{equation}
\textrm{Si}(x)=\int^{x}_{0}\frac{\sin(t)}{t}dt
\end{equation}
and $\textrm{ber}_{\nu}(x)$ is the Kelvin function defined as 
\begin{equation}
\textrm{ber}_{\nu}(x)=Re\left(J_{\nu}(xe^{3\pi i/4})\right)
\end{equation}
(here we denote 
$$
\textrm{ber}^{(1)}(x)=\frac{d}{dx}\textrm{ber}(x)=\frac{\textrm{bei}_{1}(x)}{\sqrt{2}}+\frac{\textrm{ber}_{1}(x)}{\sqrt{2}},
$$
$\textrm{ber}(x)=\textrm{ber}_{0}(x)$).\\  
Making the change of variable $tY=w^2$ we get
\begin{equation}
c(\sin^2,Y)=-\sqrt{2}\int^{\infty}_{0}\frac{\textrm{ber}^{(1)}(2\sqrt{2}t)}{e^{t^2/Y}-1}dt+\sin^2(Y)-Y\textrm{Si}(2Y)
\end{equation}
If we set $Y=Y_n=\pi(n+\theta)x$ and $Y^{*}=Y^{*}_n=\pi(n+1-\theta)x$ then using (24),(27) we get
$$
\sum^{\infty}_{m=1}\frac{J_1\left(4\pi\sqrt{m(n+\theta)x}\right)}{\sqrt{m(n+\theta)}}=\frac{\sqrt{x}}{2Y_n}+\frac{2\sqrt{x}c(\sin^2,Y_n)}{Y_n}-\frac{2\sqrt{x}}{Y_n}\int^{1}_{0}\sin^2\left(\frac{Y_n}{t}\right)dt=
$$
\begin{equation}
=\sqrt{x}\left(-\pi+\frac{1}{2Y_n}-\frac{2\sqrt{2}}{Y_n}\int^{\infty}_{0}\frac{ber^{(1)}(2\sqrt{2}t)}{e^{t^2/Y_n}-1}dt\right)
\end{equation}
Hence for $\theta=1/4$ we get
$$
P(x)=2\sqrt{x}\sum^{\infty}_{n=0}\sum^{\infty}_{m=1}\left(\frac{J_1\left(4\pi\sqrt{m(n+\frac{1}{4})x}\right)}{\sqrt{m(n+\frac{1}{4})}}-\frac{J_1\left(4\pi\sqrt{m(n+\frac{3}{4})x}\right)}{\sqrt{m(n+\frac{3}{4})}}\right)=
$$
$$
=1+4\sqrt{2}x\sum^{\infty}_{n=0}\int^{\infty}_{0}\textrm{ber}^{(1)}(2\sqrt{2}t)\left(\frac{1}{Y_n(e^{t^2/Y_n}-1)}-\frac{1}{Y^{*}_n(e^{t^2/Y^{*}_n}-1)}\right)dt
$$
Because
\begin{equation}
\sum^{\infty}_{n=0}\left(\frac{1}{2\pi (n+\theta)x^{1/2}}-\frac{1}{2\pi (n+1-\theta)x^{1/2}}\right)=\frac{\cot(\pi \theta)}{2x^{1/2}}
\end{equation}
\\
\textbf{Theorem 2.}\\
Let $Y_n=\pi(n+\theta)x$, then
\begin{equation}
\sum^{\infty}_{m=1}\frac{J_1\left(4\pi\sqrt{m(n+\theta)x}\right)}{\sqrt{m(n+\theta)}}
=\sqrt{x}\left(-\pi+\frac{1}{2Y_n}-\frac{2\sqrt{2}}{Y_n}\int^{\infty}_{0}\frac{ber^{(1)}(2\sqrt{2}t)}{e^{t^2/Y_n}-1}dt\right)
\end{equation}
\\
\textbf{Note 1.}\\ 
Extending our thoughts of Lemma 2 one can show that exists sequence $c_n$, $n=0,1,2,\ldots$ such that 
\begin{equation}
\zeta(s)-\frac{1}{s-1}=\sum^{\infty}_{n=0}\frac{c_n}{(s+1)(s+2)\ldots(s+n)}
\end{equation}
Hence
\begin{equation}
c(f,x)=c_0f(x)+\sum^{\infty}_{n=1}\frac{c_n}{x^n}\int^{x}_{0}\int^{x_{n-1}}_{0}\int^{x_{n-2}}_{0}\ldots\int^{x_1}_{0}f(t)dtdx_1\ldots dx_{n-2}dx_{n-1}
\end{equation}
And if $f(x)=\sin^2(x)$, then
$$
\frac{1}{x^n}\int^{x}_{0}\int^{x_{n-1}}_{0}\int^{x_{n-2}}_{0}\ldots\int^{x_1}_{0}f(t)dtdx_1\ldots dx_{n-2}dx_{n-1}=
$$
\begin{equation}
=\frac{1}{2\cdot n!}+P_n\left(\frac{1}{x}\right)+\frac{\{\sin^2(x),\sin(2x)\}}{x^n}
\end{equation} 
where $\textrm{deg}(P_n(x))=2n-2$ and $P(0)=0$. Clearly this series lead to a approximation method.\\
The general result is
\begin{equation}
c\left(\sin^2,x\right)=\sum^{\infty}_{k=0}\frac{(-1)^{k}A_{2k}}{2^{2k+1}x^{2k}}+\sin^2(x)\sum^{\infty}_{k=0}\frac{(-1)^kc_{2k}}{2^{2k}x^{2k}}+\frac{\sin(2x)}{2}\sum^{\infty}_{k=0}\frac{(-1)^kc_{2k+1}}{2^{2k+1}x^{2k+1}}
\end{equation}
where
\begin{equation}
A_{2k}=\sum^{\infty}_{s=1}\frac{c_{s+2k}}{s!}
\end{equation}
Set in (35) where $x$ the $Y_n=\pi\left(n+\frac{1}{4}\right)x$. It is easy to see someone that then we have $\sin^2(Y_n)=\frac{1}{2},1,\frac{1}{2},0$ and $\frac{1}{2}\sin(2Y_n)=\frac{1}{2},0,-\frac{1}{2},0$, respectively, for $x\equiv1(mod4)$, $2(mod4)$, $3(mod4)$, $0(mod4)$.\\
Assume that $x\equiv0(mod4)$, then 
\begin{equation}
c(\sin^2,Y_n)=\sum^{\infty}_{k=0}\frac{(-1)^kA_{2k}}{2^{2k+1}Y_n^{2k}}    
\end{equation}
This case shows that different modulus clases of $x$ lead to different aproximations. However we did not proceed in this way for is too dificult to evaluate the $c_n$.\\ 

Under the substitution
\begin{equation}
\zeta(2n)=\frac{(-1)^{n-1}(2\pi)^{2n}B_{2n}}{2(2n)!}\textrm{, }n=1,2,\ldots
\end{equation} 
relation (25) becomes
\begin{equation}
c(\sin^2,Y)=f(Y)+\sin^2(Y)-Y\textrm{Si}(2Y)
\end{equation}
where
\begin{equation}
f(x):=\frac{1}{4}\sum^{\infty}_{n=1}\frac{B_{2n}}{((2n)!)^2}(4\pi x)^{2n}=\sum^{\infty}_{a=1}\sin^2\left(\frac{x}{a}\right)
\end{equation}
Hence\\
\\
\textbf{Theorem 3.}\\
For $n=0,1,2,\ldots$ holds
\begin{equation}
\sum^{\infty}_{m=1}\frac{J_1\left(4\pi\sqrt{m(n+\theta)x}\right)}{\sqrt{m(n+\theta)}}=\sqrt{x}\left(-\pi+\frac{1}{2Y_n}+2\frac{f(Y_n)}{Y_n}\right)
\end{equation}
and for the error term of (10) it is
\begin{equation}
P(x)=1+4x\sum^{\infty}_{n=0}\left(\frac{f( Y_n)}{Y_n}-\frac{f(Y^{*}_n)}{Y^{*}_n}\right)
\end{equation}
where $Y_n=\pi(n+1/4)x$ and $Y^{*}_{n}=\pi(n+3/4)x$.\\
\\

Note here that $f(x)$ is absolutely convergent and the rate of convegence is about $1/a^2$.\\
Consider now the quantity
$$
S=\sum^{\infty}_{n=0}\sum^{\infty}_{m=1}\left(\frac{J_1\left(4\pi\sqrt{m(n+\theta)x}\right)}{\sqrt{m(n+\theta)}}-\frac{J_1\left(4\pi\sqrt{m(n+1-\theta)x}\right)}{\sqrt{m(n+1-\theta)}}\right)=
$$
$$
=\frac{1}{2\sqrt{x}}+2\sqrt{x}\sum^{\infty}_{n=0}\left(\frac{f(Y_n)}{Y_n}-\frac{f(Y^{*}_{n})}{Y^{*}_n}\right)=
$$
$$
=\frac{1}{2\sqrt{x}}+2\sqrt{x}\sum^{\infty}_{n=0}\left(\frac{f'(\xi_n)}{\xi_n}-\frac{f(\xi_n)}{\xi_n^2}\right)\frac{\pi x}{2}
$$
where $\xi_n\in[Y_n,Y_n^{*}]$.\\
It is easy to see someone that hold the following asymptotic expansions for the Riemann $\zeta$-zeta function.\\
\begin{equation}
\sum^{\infty}_{k=x+1}\frac{1}{k^{2n}}=\frac{1}{x^{2n-1}}+\frac{1}{(2n-1)x^{2n-1}}-\frac{1/2}{x^{2n}}+\frac{n/6}{x^{2n+1}}+O\left(x^{-2n-3}\right)\textrm{, }x\rightarrow\infty
\end{equation}
\begin{equation}
\sum^{\infty}_{k=y+1}\frac{1}{k^{2n+2}}=\frac{1}{(2n+1)y^{2n+1}}-\frac{1/2}{y^{2n+2}}+\frac{n+1}{6y^{2n+3}}+O\left(y^{-2n-4}\right)\textrm{, }y\rightarrow\infty
\end{equation}
From these we get
$$
\frac{f(x)}{x}=\frac{1}{x}\sum^{\infty}_{k=1}\sin^2\left(\frac{x}{k}\right)=\frac{1}{x}\sum^{x}_{k=1}\sin^2\left(\frac{x}{k}\right)+\frac{1}{x}\sum^{\infty}_{k=x+1}\frac{1}{2}\left(1-\cos\left(\frac{2x}{k}\right)\right)=
$$
$$
=\int^{1}_{0}\sin^2\left(\frac{x}{t}\right)dt+I_x+\sum^{\infty}_{n=1}\frac{(-1)^{n-1}2^{2n-1}x^{2n-1}}{(2n)!}\sum^{\infty}_{k=x+1}\frac{1}{k^{2n}}=
$$
$$
=\int^{1}_{0}\sin^2\left(\frac{1}{t}\right)dt+I_x+\frac{\cos(2)-1}{2}+\textrm{Si}(2)+\frac{\cos(2)-1}{4x}+\frac{\sin(2)}{12x^2}+O\left(x^{-4}\right)
$$
Hence
\begin{equation}
\frac{f(x)}{x^2}=\frac{1}{x}\int^{1}_{0}\sin^2\left(\frac{1}{t}\right)dt+\frac{I_x}{x}+\frac{\cos(2)-1+2\textrm{Si}(2)}{2x}+\frac{\cos(2)-1}{4x^2}+\frac{\sin(2)}{12x^3}+O\left(x^{-5}\right)
\end{equation}
Also
$$
\frac{f'(x)}{x}=\frac{1}{x}\sum^{\infty}_{k=1}\frac{1}{k}\sin\left(\frac{2x}{k}\right)=
$$
$$
=\frac{1}{x}\sum^{2x}_{k=1}\frac{1}{k}\sin\left(\frac{2x}{k}\right)+\frac{1}{x}\sum^{\infty}_{k=2x+1}\frac{1}{k}\sum^{\infty}_{n=0}\frac{(-1)^n2^{2n+1}x^{2n+1}}{(2n+1)!k^{2n+1}}=\ldots
$$
$$
=\frac{1}{x}\frac{1}{2x}\sum^{2x}_{k=1}\frac{2x}{k}\sin\left(\frac{2x}{k}\right)+\textrm{Si(1)}{x}-\frac{\sin(1)}{4x^2}+\frac{\cos(1)+\sin(1)}{48x^3}+O\left(x^{-5}\right)
$$
Hence
\begin{equation}
\frac{f'(x)}{x}=\frac{1}{x}\int^{1}_{0}\frac{1}{t}\sin\left(\frac{1}{t}\right)dt+\frac{II_x}{x}+\frac{\textrm{Si}(1)}{x}-\frac{\sin(1)}{4x^2}+\frac{\cos(1)+\sin(1)}{48x^3}+O\left(x^{-5}\right)
\end{equation}
Setting where $x\rightarrow\xi_n$ in (45) and (46), we can write
\begin{equation}
\frac{f'(\xi_n)}{\xi_n}-\frac{f(\xi_n)}{\xi_n^2}=\frac{1}{\xi_n}[E_1(\xi_n)-E_2(\xi_n)]+O\left(\xi^{-2}_{n}\right)
\end{equation}
where
\begin{equation}
E_1(x)=\frac{1}{x}\sum^{x}_{k=1}\sin^2\left(\frac{x}{k}\right)-\int^{1}_{0}\sin^2\left(\frac{1}{t}\right)dt
\end{equation}
and
\begin{equation}
E_2(x)=\frac{1}{2x}\sum^{2x}_{k=1}\frac{2x}{k}\sin\left(\frac{2x}{k}\right)-\int^{1}_{0}\frac{1}{t}\sin\left(\frac{1}{t}\right)dt
\end{equation}
Hence\\
\\
\textbf{Theorem 4.}
\begin{equation}
P(x)=1+2\pi x^2\sum^{\infty}_{n=0}\left(\frac{E_1(\xi_n)}{\xi_n}-\frac{E_2(\xi_n)}{\xi_n}\right)+O\left(\sum^{\infty}_{n=0}\frac{2\sqrt{x}}{\xi_n^2}\right)\textrm{, }\xi_n\in[Y_n,Y_n^{*}]
\end{equation}
\\
\textbf{Notes 2.}\\
\textbf{1.} The $E_1$ and $E_2$ are the error terms of the Riemann approximation of integrals with the usual rectangular method.\\   
\textbf{2.} It also holds the following usefull generalized expansion    
$$
\frac{1}{x}\sum^{x}_{k=1}f\left(\frac{x}{k}\right)-\int^{1}_{0}f\left(\frac{1}{t}\right)dt=-\int^{1/x}_{0}f\left(\frac{1}{t}\right)dt+\frac{f(0)}{x}+f'(0)\left(\gamma+\frac{1}{2x}-\frac{1}{12x^2}\right)+
$$
\begin{equation}
+\frac{c(f,x)}{x}+\frac{f(1)-f(0)-f'(0)}{2x}-\frac{f'(1)-f'(0)}{12x^2}+O\left(x^{-4}\right)
\end{equation}
where 
\begin{equation}
c(f,x)=\sum^{\infty}_{s=2}\frac{f^{(s)}(0)x^s}{s!}\left(\zeta(s)-\frac{1}{s-1}\right)
\end{equation}
which is generalization of Lemma 1. For to prove it one can use 
\begin{equation}
\sum^{x}_{k=1}\frac{1}{k}=\log(x)+\gamma+\frac{1}{2x}-\frac{1}{12x^2}+O\left(x^{-4}\right)\textrm{, }x\rightarrow\infty
\end{equation}
\begin{equation}
\sum^{\infty}_{k=x+1}\frac{1}{k^s}=\frac{1}{(s-1)x^{s-1}}-\frac{1}{2x^s}+\frac{s}{12x^{s+1}}+O\left(x^{-s-3}\right)\textrm{, }x\rightarrow\infty
\end{equation}
\begin{equation}
\frac{1}{x}\int^{x}_{1}f\left(\frac{1}{t}\right)dt=\frac{x-1}{x}f(0)+f'(0)\log(x)+\frac{1}{x}\sum^{\infty}_{s=2}\frac{f^{(s)}(0)x^s}{s!(s-1)}-\sum^{\infty}_{s=2}\frac{f^{(s)}(0)}{s!(s-1)}
\end{equation}
and 
\begin{equation}
\frac{1}{x}\int^{x}_{1}f\left(\frac{1}{t}\right)dt=-\int^{1/x}_{0}f\left(\frac{1}{t}\right)dt+\int^{1}_{0}f\left(\frac{1}{t}\right)dt
\end{equation}
\\

Set now $f_1(x)=\sin^2(x)$ and $f_2(x)=x\sin(x)$.
The error term $E_2(x)$ is not behave properly because has not bounded singular point at $t=0$ i.e the function $\frac{1}{t}\sin\left(\frac{1}{t}\right)$ goes to infinity when $t\rightarrow 0$. The error term $E_1(x)$ behaves nice and $E_1(x)=O\left(x^{-1}\right)\textrm{, }x\rightarrow\infty$.  
Hence if $c_2(x)=c(f_2,x)$, then 
$$
P_1=\sum^{\infty}_{n=0}\left(\frac{E_1(\xi_n)-E_2(\xi_n)}{\xi_n}\right)=\sum^{\infty}_{n=0}\frac{O(1)}{\xi^2_n}-\sum^{\infty}_{n=0}\frac{E_2(\xi_n)}{\xi_n}=O\left(\sum^{\infty}_{n=0}\frac{c(f_2,2\xi_n)}{2\xi_n^2}\right)
$$
\begin{equation}
=O\left(\sum^{\infty}_{n=0}\frac{c_2(2\xi_n)}{\xi_n^2}\right)
\end{equation}
But
$$
\frac{c_2(2x)}{2x}=\lim_{M\rightarrow\infty}\left\{\sum^{M}_{k=1}\frac{1}{k}\sin\left(\frac{2x}{k}\right)-\int^{M}_{1}\frac{1}{t}\sin\left(\frac{2x}{t}\right)dt\right\}=
$$
$$
=\sum^{\infty}_{k=1}\frac{1}{n}\sin\left(\frac{2x}{k}\right)-\int^{\infty}_{1}\frac{1}{t}\sin\left(\frac{2x}{t}\right)dt=\sum^{\infty}_{k=1}\frac{1}{k}\sin\left(\frac{2x}{k}\right)-\textrm{Si}(2x)
$$
Let $X=\lambda_n$ be the roots of the equation
\begin{equation}
\sum^{\infty}_{k=1}\frac{1}{k}\sin\left(\frac{2X}{k}\right)-\textrm{Si}(2\xi_n)=0
\end{equation}
then
$$
P_1=O\left(\sum^{\infty}_{n=0}\frac{1}{\xi_n}\left\{\sum^{\infty}_{k=1}\frac{1}{k}\sin\left(\frac{2\xi_n}{k}\right)-\textrm{Si}(2\xi_n)\right\}\right)=
$$
$$
=O\left(\sum^{\infty}_{n=0}\frac{1}{\xi_n}\left\{\sum^{\infty}_{k=1}\frac{1}{k}\sin\left(\frac{2\xi_n}{k}\right)-\sum^{\infty}_{k=1}\frac{1}{k}\sin\left(\frac{2\lambda_n}{k}\right)\right\}\right)=
$$
$$
=O\left(\sum^{\infty}_{n=0}\frac{1}{\xi_n}\left\{\sum^{\infty}_{k=1}\frac{\sin\left(\frac{2\xi_n}{k}\right)-\sin\left(\frac{2\lambda_n}{k}\right)}{k}\right\}\right)
$$
Since $|\sin(x)-\sin(y)|\leq C|x-y|$, we get
$$
P_1=O\left(\sum^{\infty}_{n=0}\sum^{\infty}_{k=1}\frac{1}{k^2}\frac{|\xi_n-\lambda_n|}{\xi_n}\right)=O\left(\sum^{\infty}_{n=0}\frac{|\xi_n-\lambda_n|}{\xi_n}\right)
$$
and holds the following conditional theorem\\
\\
\textbf{Theorem 5.}\\
If
\begin{equation}
|\xi_n-\lambda_n|=O\left(\frac{1}{\xi_n^{3/4}}\right)\textrm{, }n\rightarrow\infty
\end{equation}
then $P(x)=O\left(x^{1/4}\right)$.\\ 
\\

The sum $\sum^{\infty}_{a=1}\sin^2(x/a)$ is absolutly convergent for every $x$. Hence setting where $x\rightarrow Y_n= \pi(n+1/4)x$ and $Y^{*}_n=\pi(n+3/4)x$ and rearanging (if this is posible), $P(x)$ can be written as 
$$
P(x)=1+4x\sum^{\infty}_{k=1}\sum^{\infty}_{n=0}\left(\frac{1-\cos(\frac{2}{k}\pi(n+\frac{1}{4})x)}{2\pi(n+\frac{1}{4})x}-\frac{1-\cos(\frac{2}{k}\pi(n+\frac{3}{4})x)}{2\pi(n+\frac{3}{4})x}\right),
$$
which is not convergent since
\begin{equation}
\sum^{\infty}_{n=0}\frac{\cos(\pi(n+\frac{1}{4})x)}{\pi(n+\frac{1}{4})}=Re\left(\phi_1\left(e^{-i\pi x/4}\right)\right)
\end{equation}
where
\begin{equation}
\phi_1(x)=\frac{2}{\pi}(\arctan(x)+\textrm{arctanh}(x))
\end{equation}
and
\begin{equation}
\sum^{\infty}_{n=0}\frac{\cos(\pi(n+\frac{3}{4})x)}{\pi(n+\frac{3}{4})}=Re\left(\phi_2\left(e^{-i\pi x/4}\right)\right)
\end{equation}
where 
\begin{equation}
\phi_2(x)=\frac{2}{\pi}(-\arctan(x)+\textrm{arctanh}(x))
\end{equation} 
$$
P(x)=1+2 \sum^{\infty}_{k=1}\left\{1-Re\left(\phi_1\left(e^{- \frac{i\pi x}{2k}}\right)\right)-Re\left(\phi_2\left(e^{-\frac{i\pi x}{2k}}\right)\right)\right\}
$$
Hence
\begin{equation}
P(x)=1+2 Re\sum^{\infty}_{k=1}\left(1-\frac{4}{\pi }\arctan\left(e^{-\frac{i\pi x}{2k}}\right)\right),
\end{equation}
which from my point of view seems to have no meaning since it is 1.\\
\\

Assume now that we wish to approximate (25). By taking finite number of terms in (34) (here 3) we get a bounded error term up to a constant due to $\frac{1}{2\cdot n!}$ in (35). But if instead consider the difference in (10), then with the approximation method (33)-(34)-(35), the three first terms give the maximum bound
$$
S=2\sqrt{x}\sum^{\infty}_{n=0}\sum^{\infty}_{m=1}\left(\frac{J_1(4\pi\sqrt{m(n+\frac{1}{4})x})}{\sqrt{m(n+\frac{1}{4})}}-\frac{J_1(4\pi\sqrt{m(n+\frac{3}{4})x})}{\sqrt{m(n+\frac{3}{4})}}\right)=
$$
$$
=-\sum^{\infty}_{n=0}\left(\frac{\sqrt{x}}{16Y_n}-\frac{\sqrt{x}}{16Y^{*}_n}\right)-
$$
$$
-\sum^{\infty}_{n=0}\frac{\sqrt{x}\cos(2Y_n)}{2Y_n}+\sum^{\infty}_{n=0}\frac{\sqrt{x}\cos(2Y^{*}_n)}{2Y^{*}_n}+O(x^{-2})=O(1)
$$
For which
\begin{equation}
\sum^{\infty}_{n=0}\frac{\cos(2\pi(n+\frac{1}{4})x)}{\pi(n+\frac{1}{4})x}-\sum^{\infty}_{n=0}\frac{\cos(2\pi(n+\frac{3}{4})x)}{\pi(n+\frac{3}{4})x}=O(1)
\end{equation}
and according the cutting of the series, the error is   
\begin{equation}
O\left(\frac{1}{Y_n^4}\int^{Y_n}_{0}\int^{w}_{0}\int^{z}_{0}\sin^2(t)dtdzdw-\frac{1}{Y^{{*}4}_n}\int^{Y^{*}_n}_{0}\int^{w}_{0}\int^{z}_{0}\sin^2(t)dtdzdw\right)=O(1),
\end{equation}
since $f(x)=\sin^2(x)$,  
$$
\textrm{Si}(x)\sim\frac{\pi}{2}-\frac{\sin(x)}{x}\left(\sum^{\infty}_{n=0}\frac{(-1)^n(2n+1)!}{x^{2n+1}}\right)-\frac{\cos(x)}{x}\left(\sum^{\infty}_{n=0}\frac{(-1)^n(2n)!}{x^{2n}}\right)
$$
as $x\rightarrow\infty$.\\
If we assume that the method of expansion (33) leads to approximation of $c(f,x)$ then from the above we must have
\begin{equation}
\sum^{*}_{0\leq n\leq x}r_2(n)=\pi x+O(1)\textrm{, }x\rightarrow\infty 
\end{equation}    
Note that we use $c_0=1,c_1=-1,c_2=-1/4$ and
$$
c(f,x)=f(x)-\frac{1}{x}\int^{x}_{0}f(t)dt-\frac{1/4}{x^2}\int^{x}_{0}\int^{w}_{0}f(t)dtdw+
$$
$$
+O\left(\frac{1}{x^4}\int^{x}_{0}\int^{w}_{0}\int^{z}_{0}f(t)dtdzdw\right)
$$
and
$$
S=\sum^{\infty}_{n=0}[\left\{\frac{\sqrt{x}}{2Y_n}+\frac{\sqrt{x}}{Y_n}c(f,Y_n)-\frac{\sqrt{x}}{Y_n}\int^{1}_{0}f\left(\frac{Y_n}{t}\right)dt\right\}-
$$
$$
-\left\{\frac{\sqrt{x}}{2Y^{*}_n}+\frac{\sqrt{x}}{Y^{*}_n}c(f,Y^{*}_n)-\frac{\sqrt{x}}{Y^{*}_n}\int^{1}_{0}f\left(\frac{Y^{*}_n}{t}\right)dt\right\}]
$$
with $f(x)=\sin^2(x)$.
\[
\]

\centerline{\bf References}\vskip .2in

\noindent

[1]: G.E. Andrews, B.C. Berndt. 'Ramanujan's Lost Notebook Part IV'. Springer., New York, Heidelberg, Dordrecht, London. 2013.

[2]: Nikos Bagis, M.L. Glasser. 'Integrals and Series Resulting from Two Sampling Theorems'. Sampling Theory in Singnal and Image Processing., Sampling Publishing, Vol. 5, No. 1, 2006.

[3]: M. Abramowitz and I.A. Stegun. 'Handbook of Mathematical Functions'. Dover Publications, New York. 1972.
\end{document}